\documentclass{ifacconf}

\usepackage{natbib} 
\usepackage[latin1]{inputenc}
\usepackage[american]{babel}
\usepackage{graphicx}
\usepackage{amscd}
\usepackage{amsfonts}
\usepackage{amsmath}
\usepackage{amssymb}
\usepackage{helvet}
\usepackage{rotating}
\usepackage{epsfig}
\usepackage{color}
\usepackage[european]{circuitikz}
\usepackage{longtable}
\usepackage{algorithm}
\usepackage{algorithmic}
\newtheorem{theorem}{Theorem}
\newtheorem{proposition}[theorem]{Proposition}

\newtheorem{conclusion}[theorem]{Conclusion}
\newtheorem{definition}[theorem]{Definition}
\newtheorem{assumption}[theorem]{Assumption}
\newtheorem{remark}[theorem]{Remark}

\newcommand{\R}{\mathbb{R}}
\newcommand{\N}{\mathbb{N}}

\newcommand{\cK}{\mathcal{K}}
\newcommand{\cKL}{\mathcal{KL}}
\newcommand{\cB}{\mathcal{B}}

\newcommand{\horizon}{N}
\newcommand{\state}{x}
\newcommand{\statestar}{\state^\star}
\newcommand{\stateequilibrium}{\state_{*}}
\newcommand{\statedisturbed}{\overline{\state}}
\newcommand{\statezero}{\state_0}
\newcommand{\statezerodisturbed}{\overline{\state}_0}
\newcommand{\statedisturbance}{\delta_{\state}}
\newcommand{\statedisturbedmax}{\Delta_{\state}}
\newcommand{\control}{u}
\newcommand{\controlstar}{\control^\star}
\newcommand{\controlequilibrium}{\control_{*}}
\newcommand{\controldisturbed}{\overline{\control}}

\newcommand{\parameter}{p}
\newcommand{\parameterstar}{\parameter^\star}
\newcommand{\parameterequilibrium}{\parameter_{*}}
\newcommand{\parameterdisturbed}{\overline{\parameter}}
\newcommand{\parameterdisturbedmax}{\Delta_{\parameter}}
\newcommand{\parameterdisturbance}{\delta_{\parameter}}

\newcommand{\feedback}{\mu}
\newcommand{\feedbackdisturbed}{\overline{\feedback}}
\newcommand{\feedbackN}{\feedback_{\horizon}}
\newcommand{\feedbackdisturbedN}{\feedbackdisturbed_{\horizon}}

\newcommand{\distance}{d}
\newcommand{\dynamic}{f}
\newcommand{\dynamicdisturbed}{\overline{\dynamic}}

\newcommand{\dynamicdisturbance}{\delta_{\dynamic}}
\newcommand{\stagecost}{\ell}

\newcommand{\costfunction}{J}
\newcommand{\costfunctionN}{\costfunction_{\horizon}}
\newcommand{\costfunctioninfty}{\costfunction_\infty}
\newcommand{\valuefunction}{V}
\newcommand{\valuefunctionN}{\valuefunction_{\horizon}}

\newcommand{\valuefunctioninfty}{\valuefunction_\infty}

\newcommand{\valuefunctioninftyfeedbackdisturbedN}{\valuefunctioninfty^{\feedbackdisturbedN}}
\newcommand{\valuefunctioninftydisturbedfeedbackdisturbedN}{\overline{\valuefunction}_\infty^{\feedbackdisturbedN}}
\newcommand{\lipschitzconstant}{L}
\newcommand{\lipschitzconstantvaluefunction}{\lipschitzconstant_{\valuefunction}}
\newcommand{\lipschitzconstantstagecost}{\lipschitzconstant_{\stagecost}}
\newcommand{\lipschitzconstantcontrol}{\lipschitzconstant_{\control}}
\newcommand{\lipschitzconstantstate}{\lipschitzconstant_{\state}}

\newcommand{\jumpconstant}{S}
\newcommand{\jumpconstantvaluefunction}{\jumpconstant_{\valuefunction}}
\newcommand{\jumpconstantstagecost}{\jumpconstant_{\stagecost}}
\newcommand{\closedloopindex}{n}
\newcommand{\openloopindex}{k}
\newcommand{\stateset}{X}
\newcommand{\controlset}{U}
\newcommand{\controlsetopenloop}{\controlset^\indexset}
\newcommand{\parameterset}{P}
\newcommand{\parametersetopenloop}{\parameterset^\indexset}
\newcommand{\statesetconstrained}{\mathbb{\stateset}}
\newcommand{\controlsetconstrained}{\mathbb{\controlset}}
\newcommand{\indexset}{\mathcal{I}}
\newcommand{\admissiblecontrolset}{\mathcal{\controlset}^\indexset}
\newcommand{\updateset}{A}
\newcommand{\feasibilityset}{\mathcal{F}}
\newcommand{\practicalset}{\mathcal{P}}
\newcommand{\unknownset}{\mathcal{L}}

\newcommand{\realtimeerrorconstant}{c_1}
\newcommand{\realtimeerrorlinear}{c_2}

\definecolor{juergen}{rgb}{.1,.4,.9}
\definecolor{johannes}{rgb}{.9,.1,.4}

\begin{document}
\begin{frontmatter}

\title{A General Framework for Nonlinear Model Predictive Control with Abstract Updates} 
\thanks[footnoteinfo]{A basic version of this paper was presented at the IFAC NMPC'12 meeting. Corresponding author J.~Pannek. Tel. +49-89-6004-2567. Fax +49-89-6004-2136.}
\author[UniBW]{J.~Pannek}
\author[UniBW]{J.~Michael} 
\author[UniBW]{M.~Gerdts}
\address[UniBW]{Faculty of Aerospace Engineering, University of the Federal Armed Forces, 85577 Munich, Germany (johannes.michael@unibw.de, matthias.gerdts@unibw.de, juergen.pannek@unibw.de).}

\begin{abstract}
: Considering nonlinear processes which are subject to unknown but measurable disturbances, we provide both stability and feasibility proofs for nonlinear model predictive controllers with abstract updates without the use of stabilizing terminal constraints or cost. The result utilizes a relaxed Lyapunov inequality for the nominal system and reasonable affine Lipschitz conditions on the open loop change of the optimal value function and the stage costs. For this methodology, we provide proofs that the known MPC updating techniques based on sensitivities, realtime iterations and hierachically structured MPC updates satisfy our assumptions revealing a common stability framework. To illustrate our approach we present a quartercar example and show performance improvements of the updated solution with respect to comfort and handling properties.
\end{abstract}


\end{frontmatter}

\section{Introduction}
Within the last decades, model predictive control (MPC) has grown mature for both linear and nonlinear systems, see, e.g., \cite{BitmeadGeversWertz1990, CamachoBordons2004} and \cite{GruenePannek2011, RawlingsMayne2009}. Although analytically and numerically challenging, this method of approximating an infinite horizon optimal control is attractive due to its simplicity: in each sampling interval, based on a measurement of the current state, a truncated finite horizon optimal control is computed and the first element (or sometimes also more) of the resulting optimal control sequence is applied to the process as input for the next sampling interval(s). Then, the entire problem is shifted forward in time rendering the scheme to be iteratively applicable. \\
As a consequence of the truncation of the infinite horizon, stability and optimality of the closed loop may be lost. Yet, stability of the resulting closed loop can be guaranteed, either by imposing terminal point constraints as shown in \cite{Alamir2006, KeerthiGilbert1988} or Lyapunov type terminal costs and terminal regions, see \cite{ChenAllgoewer1998, MayneRawlingsRaoScokaert2000}. Apart from these modifications, an alternative approach shown in \cite{GrueneRantzer2008} utilizes a relaxed Lyapunov condition which can be shown to hold if the system is controllable in terms of the stage costs, cf. \cite{Gruene2009}. Additionally, this method allows for computing an estimate on the degree of suboptimality with respect to the infinite horizon controller, see also \cite{NevisticPrimbs1997, ShammaXiong1997} for earlier works on this topic and \cite{GruenePannek2009} for a method to compute this degree.\\
In this work, we extend the third approach to the case of parametric control systems subject to measured disturbances including subsequent disturbance rejection updates. In order to be as general as possible, an abstract update law is considered to handle the disturbances and a respective algorithm is presented. While feasibility of the resulting closed loop can be shown via self--concordant set conditions, stability requires affine Lipschitz conditions on the open loop change of the optimal value function as well as the stage cost function. Combining these conditions with the relaxed Lyapunov condition for the nominal open loop solution allows us to show a generalization of the stability proof given in \cite{PannekGerdts2012b} for the closed loop solution despite the presence of disturbances.\\
Bringing the general framework to life, we show that the known MPC update laws utilizing realtime iterations \cite{DiehlBockSchloeder2005} as well as hierarchical approaches \cite{BockDiehlKostinaSchloeder2007} and sensitivity information \cite{ZavalaBiegler2009} are special cases of our abstract update law. Despite the fact that these methods where designed for MPC problems with Lyapunov type terminal costs and terminal regions to ensure stability of the closed loop, these modifications are not required to prove that the stated affine Lipschitz conditions hold and the closed loop to be stable.\\
The paper is organized as follows: The problem formulation and the concept of practical stability are defined in Section \ref{Section:Setup and Preliminaries}. Additionally, the basic MPC algorithm with abstract updates is introduced. In the subsequent Section \ref{Section:Stability and Feasibility}, conditions for feasibility and practical stability of the closed loop are given and respective proofs are shown. Thereafter, the stability conditions are verified for the known sensitivity based update law, realtime iteration and the hierarchical MPC approach in Section \ref{Section:Update Techniques}. Moreover, we present simulation results for a quartercar subject to disturbed measurements of both the state and the sensor inputs in Section \ref{Section:Numerical Results}. Concluding our paper, we shortly summarize the results and point out areas of future research.

\section{Setup and Preliminaries}
\label{Section:Setup and Preliminaries}

The set $\N$ denotes the natural numbers, $\N_0$ the natural numbers including zero, $\R$ denotes the real numbers and $\R_0^+$ the nonnegative reals. We call a continuous function $\rho: \R_0^+ \rightarrow \R_0^+$ a class $\cK_\infty$-function if it satisfies $\rho(0) = 0$, is strictly increasing and unbounded. A continuous function $\beta: \R_0^+ \times \R_0^+ \rightarrow \R_0^+$ is said to be of class $\cKL$ if for each $r > 0$ the limit $\lim_{t \rightarrow \infty} \beta(r, t) = 0$ holds and for each $t \geq 0$ the condition $\beta(\cdot, t) \in \cK_\infty$ is satisfied. 

Throughout this paper we consider plant models driven by the dynamics
\begin{align}
	\label{Setup and Preliminaries:eq:control system}
	\state(\closedloopindex+1) = \dynamic(\state(\closedloopindex),\control(\closedloopindex), \parameter(\closedloopindex))
\end{align}
where $\state$ denotes the state of the system, $\control$ the external control and $\parameter$ a parameter which can be measured. These variables are elements of the respective state, control and parameter spaces $(\stateset, \distance_{\stateset})$, $(\controlset, \distance_{\controlset})$ and $(\parameterset, \distance_{\parameterset})$. Here, we allow $\stateset$, $\controlset$ and $\parameterset$ to be arbitrary metric spaces which renders our results applicable to discrete time dynamics induced by a sampled finite or infinite dimensional system. In this setting, constraints on both state and control are introduced by considering suitable subsets $\statesetconstrained \subseteq \stateset$ and $\controlsetconstrained \subseteq \controlset$. \\
In most cases, simulated and true trajectories do not coincide. Here, we denote the \textit{simulated} or \textit{nominal} state trajectory of \eqref{Setup and Preliminaries:eq:control system} corresponding to an initial state $\statezero \in \stateset$, a control sequence $\control = \left( \control(\openloopindex) \right)_{\openloopindex \in \N_0}$ and a nominal parameter sequence $\parameter = \left( \parameter(\openloopindex) \right)_{\openloopindex \in \N_0}$ by $\state(\cdot) = \state_{\control, \parameter}(\cdot; \statezero)$.  The true system dynamics, however, are given by
\begin{align}
	\label{Setup and Preliminaries:eq:plant}
	\statedisturbed(\closedloopindex+1) & = \dynamicdisturbed(\statedisturbed(\closedloopindex),\controldisturbed(\closedloopindex), \parameterdisturbed(\closedloopindex)) \\
	& = \dynamic(\state(\closedloopindex) + \statedisturbance(\closedloopindex),\controldisturbed(\closedloopindex), \parameter(\closedloopindex) + \parameterdisturbance(\closedloopindex)) + \dynamicdisturbance(\closedloopindex) \nonumber
\end{align}
where $\statedisturbance$, $\parameterdisturbance$ correspond to measurement errors and $\dynamicdisturbance$ represents the plant--model mismatch. To distinguish between solutions resulting from \eqref{Setup and Preliminaries:eq:control system} and \eqref{Setup and Preliminaries:eq:plant}, $\statedisturbed(\cdot) = \statedisturbed_{\controldisturbed, \parameterdisturbed}(\cdot; \statezerodisturbed)$ denotes the \textit{true} or \textit{disturbed} solution of \eqref{Setup and Preliminaries:eq:plant} which is subject to a disturbed initial value $\statezerodisturbed$, a possibly modified control sequence $\controldisturbed = \left( \controldisturbed(\openloopindex) \right)_{\openloopindex \in \N_0}$ and the disturbed parameter sequence $\parameterdisturbed = \left( \parameterdisturbed(\openloopindex) \right)_{\openloopindex \in \N_0}$. We call a point $\stateequilibrium \in \stateset$ a controlled equilibrium of \eqref{Setup and Preliminaries:eq:plant} if there exist a control $\controlequilibrium \in \controlset$ and a nominal parameter $\parameterequilibrium \in \parameterset$ satisfying $\dynamicdisturbed(\stateequilibrium, \controlequilibrium, \parameterequilibrium) = \stateequilibrium$. \\
To shorten notation, we define the abbreviation $\| x \|_{y} := \distance_i(x, y)$ for $\distance_i \in \{ \distance_\stateset, \distance_\controlset, \distance_\parameterset \}$. This allows us to formulate the following assumption which we suppose to hold throughout the paper:

\begin{assumption}\label{Setup and Preliminaries:ass:bounded disturbance}
	There exist a control forward invariant set $\feasibilityset \subset \statesetconstrained$ and constants $\statedisturbedmax, \parameterdisturbedmax > 0$ 
	such that for all $\state \in \feasibilityset$ we have $\| \statedisturbance \| = \| \statedisturbed \|_{\state} \leq \statedisturbedmax$ and additionally the bounds $\| \parameterdisturbance \| = \| \parameterdisturbed \|_{\parameter} \leq \parameterdisturbedmax$ and $\| \parameter(\closedloopindex+1) \|_{\parameter(\closedloopindex)} \leq \parameterdisturbedmax$ hold for all $\closedloopindex \in \N_0$.
\end{assumption}

For this setting, we wish to design a feedback controller based on the model \eqref{Setup and Preliminaries:eq:control system} which despite the measurement and modelling errors semiglobally practically asymptotically stabilizes the plant \eqref{Setup and Preliminaries:eq:plant} at a controlled equilibrium, i.e. which satisfies the following ISS property:

\begin{definition}\label{Setup and Preliminaries:def:stability}
	Suppose $\stateequilibrium \in \statesetconstrained$ is a controlled equilibrium of \eqref{Setup and Preliminaries:eq:plant}. Then a control sequence $\control = (\control(\closedloopindex))_{\closedloopindex \in \N_0}$ {\em semiglobally practically asymptotically stablizes} system \eqref{Setup and Preliminaries:eq:plant} if there exist $r>0$, $\beta \in \cKL$ and $\gamma \in \cK_\infty$ such that the inequality
	\begin{align}
		\label{Setup and Preliminaries:def:stability:eq1}
		\| \statedisturbed_{\control, \parameterdisturbed}(\closedloopindex, \statezerodisturbed) \|_{\stateequilibrium} \leq \beta( \| \statezerodisturbed \|_{\stateequilibrium}, \closedloopindex) + \gamma(\max\{\statedisturbance, \parameterdisturbance, \dynamicdisturbance\})
  	\end{align}
	holds for all $\statezero \in \cB_r(\stateequilibrium)$, all $\closedloopindex \in \N_0$ and all $\statedisturbance$, $\parameterdisturbance$, $\dynamicdisturbance$ satisfying Assumption \ref{Setup and Preliminaries:ass:bounded disturbance}. 
\end{definition}

To accomplish this task, we propose a two stage feedback. In the first stage, we apply model predictive control (MPC) to generate an approximation of a minimizer of the infinite horizon cost functional $\costfunctioninfty ( \state, \control, \parameter ) = \sum_{\closedloopindex=0}^\infty \stagecost(\state(\closedloopindex), \control(\closedloopindex), \parameter(\closedloopindex))$. As indicated by the notation, the basis of this computation is the model \eqref{Setup and Preliminaries:eq:control system} together with a nominal parameter sequence $\left( \parameter(\closedloopindex) \right)_{\closedloopindex \in \N_0}$. This step is applied in an advanced step setting, cf. \cite{FindeisenAllgoewer2004, ZavalaBiegler2009}, resulting in a control that is to be implemented at some future time instant. The second stage utilizes the precomputed control to locally update it based on newly obtained measurements of $\state$ and $\parameter$. Examples of updating techniques include sensitivities which have been analysed in the nonlinear optimization context for quite some time, see, e.g., \cite{GroetschelKrumkeRambau2001}, but have also been applied in the MPC context, cf. \cite{ZavalaBiegler2009}, as well as realtime iterations \cite{DiehlBockSchloeder2005}, and hierarchical MPC \cite{BockDiehlKostinaSchloeder2007}.\\
Within the first stage, we assume the stage cost $\stagecost: \stateset \times \controlset \times \parameterset \to \R_0^+$ to be continuous and to satisfy $\stagecost(\stateequilibrium, \controlequilibrium, \parameterequilibrium) = 0$ and $\stagecost(\state, \control, \parameter) > 0$ for all $\control \in \controlset$ for each $\state \not = \stateequilibrium$ and each $\parameter \not = \parameterequilibrium$. Due to the possible existence of constraints, we only consider the set of admissible controls within the minimization. To formally define this set, we introduce $\controlsetopenloop$ and $\parametersetopenloop$ denoting the set of all control and parameter sequences where $\indexset := \{ 0, 1, \ldots, \horizon - 1\}$ and $\horizon \in \N \cup \{ \infty \}$ specifies the length of these sequences. Then, the set of admissible controls for a fixed parameter sequence $\parameter \in \parametersetopenloop$ is given by
\begin{align*}
	\admissiblecontrolset(\statezero, \parameter) := \{ \control \in \controlsetopenloop \mid & \dynamic (\state(\closedloopindex), \control(\closedloopindex), \parameter(\closedloopindex)) \in \statesetconstrained \; \text{and} \; \\
	& \control(\closedloopindex) \in \controlsetconstrained \; \text{for all} \; \closedloopindex \in \indexset \}.
\end{align*}
Since solving the infinite horizon optimal control problem is in most cases computationally intractable, the idea of MPC is to approximate such a solution via a series of finite horizon problems. To this end, the truncated cost functional
\begin{align}
	\label{Setup and Preliminaries:eq:cost functional}
	\costfunctionN(\state, \control, \parameter) := \sum_{\openloopindex=0}^{\horizon-1} \stagecost(\state(\openloopindex), \control(\openloopindex), \parameter(\openloopindex)) 
\end{align}
with horizon length $N \in \N_{\geq 2}$ is minimized revealing a finite optimal control sequence $\controlstar(\cdot, \state, \parameter)$. Then, only the first element is implemented defining the feedback law $\feedbackN(\state, \parameter) := \controlstar(0, \state, \parameter)$ and the optimization horizon is shifted forward in time which allows the method to be iteratively applicable, see also \cite{Alamir2006, CamachoBordons2004, GruenePannek2011, RawlingsMayne2009} for further details. Incorporating the second stage of the proposed feedback design, the only modification of the MPC methodology is to compute an updated control $\feedbackdisturbedN(\statedisturbed, \parameterdisturbed)$ to replace $\feedbackN(\state, \parameter)$ upon implementation. Here, we assume such an update to be instantly computable. \\
As a result we obtain the following Algorithm \ref{Setup and Preliminaries:alg:mpc} where we used the short notation $\parameter_{\closedloopindex} := (\parameter(\closedloopindex+\openloopindex))_{\openloopindex \in \indexset}$ abbreviating the parameter sequence used for computing or updating the open loop control at time instant $\closedloopindex$.

\begin{algorithm}
	\caption{MPC Algorithm}\label{Setup and Preliminaries:alg:mpc}
	\begin{enumerate}
		\item Obtain measurement of $\statedisturbed(\closedloopindex)$ and $\parameterdisturbed_{\closedloopindex}$
		\item Update control $\controlstar(\cdot, \state(\closedloopindex), \parameter_{\closedloopindex})$ to $\controldisturbed(\cdot, \statedisturbed(\closedloopindex), \parameterdisturbed_{\closedloopindex})$ and apply control $\feedbackdisturbedN(\statedisturbed(\closedloopindex), \parameterdisturbed_{\closedloopindex}) := \controldisturbed(0, \statedisturbed(\closedloopindex), \parameterdisturbed_{\closedloopindex})$
		\item Set $\state(\closedloopindex) := \statedisturbed(\closedloopindex)$ and $\parameter_{\closedloopindex+1} := \parameterdisturbed_{\closedloopindex}$, predict $\state(\closedloopindex+1)$ using dynamic \eqref{Setup and Preliminaries:eq:control system} together with $\statedisturbed(\closedloopindex)$, $\feedbackdisturbedN(\statedisturbed(\closedloopindex), \parameterdisturbed_{\closedloopindex})$ and $\parameterdisturbed(\closedloopindex)$
		\item Compute $\controlstar(\cdot, \state(\closedloopindex+1), \parameter_{\closedloopindex+1})$, set $\closedloopindex := \closedloopindex + 1$ and goto step (1)
	\end{enumerate}
\end{algorithm}

\begin{remark}
	Within Algorithm \ref{Setup and Preliminaries:alg:mpc} a sequence of nominal future parameters $\parameter$ is necessary to evaluate both the dynamic and the costfunction. A respective sequence can be computed by extrapolation methods or measured via forward sensors such as road scanning laser sensors in a car or thermometers at intake pipes of a chemical plant.
\end{remark}

Here, given $\horizon \in \N$ we suppose that for each $\state \in \statesetconstrained$ and $\parameter \in \parametersetopenloop$ there exists a minimizer $\controlstar(\cdot, \state, \parameter) \in \admissiblecontrolset$ of \eqref{Setup and Preliminaries:eq:cost functional}. Note that this assumption guarantees the existence of a feasible solution for each $\state \in \statesetconstrained$ and each $\parameter \in \parametersetopenloop$, see also \cite{KerriganMaciejowski2001, PrimbsNevistic2000} for relaxed conditions. Hence, we can define the optimal value function on a finite horizon via
\begin{align}
	\label{Setup and Preliminaries:eq:value function}
	\valuefunctionN(\state, \parameter) := \min_{\control \in \admissiblecontrolset(\state, \parameter)} \costfunctionN(\state, \control, \parameter).
\end{align}
If additionally there exists a feasible updated MPC feedback $\feedbackdisturbedN(\statedisturbed, \parameterdisturbed)$ then we define the associated closed loop costs associated via
\begin{align}
	\label{Setup and Preliminaries:eq:value function disturbed}
	\valuefunctioninftyfeedbackdisturbedN(\statedisturbed, \parameterdisturbed)  := \sum_{\closedloopindex=0}^{\infty} \stagecost(\statedisturbed(\closedloopindex), \feedbackdisturbedN(\statedisturbed(\closedloopindex), \parameterdisturbed_{\closedloopindex}), \parameterdisturbed(\closedloopindex)).
\end{align}
Instead of showing semiglobal practical asymptotic stability of system \eqref{Setup and Preliminaries:eq:plant} via ISS Lyapunov functions, cf. \cite{JiangWang2001, MagniScattolini2007}, we utilize the concept of $\practicalset$--practical asymptotic stability and ``truncated'' Lyapunov functions. As outlined in \cite[Chapter 8.5]{GruenePannek2011}, ISS of a system \eqref{Setup and Preliminaries:eq:plant} can be shown via a suitable choice of $\gamma \in \cK_\infty$ in \eqref{Setup and Preliminaries:def:stability:eq1} if the following holds:

\begin{definition}\label{Setup and Preliminaries:def:p practical}
	Let $\updateset \subset \statesetconstrained$ be a forward invariant set with respect to all possible disturbances satisfying Assumption \ref{Setup and Preliminaries:ass:bounded disturbance} and let $\practicalset \subset \updateset$. Then $\statestar \in \practicalset$ is $\practicalset$--practically asymptotically stable on $\updateset$ if there exists $\beta \in \cKL$ such that
	\begin{align}
		\label{Setup and Preliminaries:eq:closed loop}
		\| \statedisturbed(\closedloopindex) \|_{\statestar} \leq \beta( \| \statezerodisturbed \|_{\statestar}, \closedloopindex)
	\end{align}
	holds for all $\statezerodisturbed \in \updateset$ and all $n \in \N_0$ with $\statedisturbed(\closedloopindex) \not \in \practicalset$.
\end{definition}

If we limit the Lyapunov property to hold outside the practical region $\practicalset$ and define $\| \parameterdisturbed \|_{\parameter}^{\max} := \max_{\closedloopindex \in \N_0} \| \parameterdisturbed(\closedloopindex) \|_{\parameter(\closedloopindex)}$, then $\practicalset$--practical asymptotic stability can be guaranteed, cf. \cite[Theorem 2.20]{GruenePannek2011} for a corresponding proof.

\begin{theorem}\label{Setup and Preliminaries:thm:p practical}
	Suppose  $\updateset \subset \statesetconstrained$ is a forward invariant set with respect to disturbances satisfying Assumption \ref{Setup and Preliminaries:ass:bounded disturbance}, $\practicalset \subset \updateset$ and $\statestar \in \practicalset$. If there exist $\alpha_1$, $\alpha_2$, $\alpha_3 \in \cK_\infty$ and a Lyapunov function $\valuefunction$ on $D = \updateset \setminus \practicalset$ satisfying
	\begin{align*}
		\alpha_1( \| \statedisturbed \|_{\statestar} + \| \parameterdisturbed \|_{\parameterstar}^{\max} ) \leq \valuefunction(\statedisturbed, \parameterdisturbed) \leq  \alpha_2( \| \statedisturbed \|_{\statestar} + \| \parameterdisturbed \|_{\parameterstar}^{\max} ) \\
		\valuefunction(\statedisturbed, \parameterdisturbed) \geq \valuefunction(\dynamic(\statedisturbed, \feedbackdisturbedN(\statedisturbed, \parameterdisturbed), \parameterdisturbed), \parameterdisturbed) - \alpha_3( \| \statedisturbed \|_{\statestar} )
	\end{align*}
	then $\statestar$ is $\practicalset$--practically asymptotically stable on $\updateset$.
\end{theorem}

\section{Stability and Feasibility}
\label{Section:Stability and Feasibility}

Before approaching the stability problem of the proposed Algorithm \ref{Setup and Preliminaries:alg:mpc}, we first need to guarantee existence of a feasible updated control $\feedbackdisturbedN(\statedisturbed, \parameterdisturbed)$.

\begin{theorem}\label{Stability and Feasibility:thm:feasibility}
	Suppose Assumption \ref{Setup and Preliminaries:ass:bounded disturbance} holds. If there exist a set $\updateset \subset \statesetconstrained$ and a feedback $\feedbackN(\state, \parameter)$ such that
	\begin{itemize}
		\item[(i)] $\inf_{x \in \feasibilityset, y \in \stateset \setminus \updateset} \distance_{\stateset}(x, y) \geq \statedisturbedmax$ holds and
		\item[(ii)] for each $\state \in \feasibilityset$ we have $\dynamic(\state,\feedbackN(\state, \parameter), \parameter) \in \feasibilityset$,
	\end{itemize}
	then the closed loop solution resulting from Algorithm \ref{Setup and Preliminaries:alg:mpc} is feasible for each initial value $\statezero \in \feasibilityset$, i.e. $\statedisturbed(\closedloopindex) \in \statesetconstrained$ for all $\closedloopindex \in \N_0$.
\end{theorem}
\textbf{Proof:} To prove the assertion by induction, first we consider estimates $\state(0) = \statezero \in \feasibilityset$ of the initial value and $\parameter_0 \in \parameterset$ of the parameter sequence to be given. Then, by Assumption \ref{Setup and Preliminaries:ass:bounded disturbance} we have $\| \statedisturbed(0) \|_{\state(0)} \leq \statedisturbedmax$ and by (i) it follows that $\statedisturbed(0) \in \updateset \subset \statesetconstrained$. In the induction step we consider estimates $\state(\closedloopindex) \in \feasibilityset$ and $\parameter_{\closedloopindex} \in \parameterset$. Now, we can use (ii) to obtain $\state(\closedloopindex+1) = \dynamic(\state(\closedloopindex),\feedbackN(\state(\closedloopindex), \parameter_{\closedloopindex}), \parameter_{\closedloopindex}) \in \feasibilityset$. Hence, Assumption \ref{Setup and Preliminaries:ass:bounded disturbance} reveals we have $\| \statedisturbed(\closedloopindex+1) \|_{\state(\closedloopindex+1)} \leq \statedisturbedmax$ and again (i) allows us to conclude $\statedisturbed(\closedloopindex) \in \updateset \subset \statesetconstrained$. \qed

\begin{remark}
	We like to note that Assumption \ref{Setup and Preliminaries:ass:bounded disturbance} is restrictive in terms of one step boundedness of the dynamics. While this condition may be weakened easily, it does not make sense from an updating technique point of view for the control which requires a certain closeness of nominal and disturbed solution.
\end{remark}

Now that we have certified existence of a solution resulting from Algorithm \ref{Setup and Preliminaries:alg:mpc} we show that the obtained feedback semiglobally practically asymptotically stabilizes system \eqref{Setup and Preliminaries:eq:plant} in the sense of Definition \ref{Setup and Preliminaries:def:stability}. In contrast to approaches based on terminal constraints or Lyapunov type terminal costs, cf., e.g., \cite{DiehlBockSchloeder2005, BockDiehlKostinaSchloeder2007, ZavalaBiegler2009}, we consider a relaxed Lyapunov condition to hold for the nominal case, cf. \cite{GrueneRantzer2008, LincolnRantzer2006}. This property is not artificial but can be shown to hold if $\horizon$ is chosen sufficiently large, see \cite{AlamirBornard1995, GrimmMessinaTeel2005, JadbabaieHauser2005, PannekWorthmann2010}.\\
One way to prove a similar result in the disturbed case of system \eqref{Setup and Preliminaries:eq:plant} is to modify the stage cost $\stagecost$ to be positive definite with respect to a robustly stabilizable forward invariant neighbourhood of $\statestar$. The computation of such a neighbourhood, however, may be impossible. Here, we consider the stage cost $\stagecost$ to be positive definite with respect to $\statestar$ only, i.e. we ignore the effects of disturbances. Since the stage cost typically decrease towards the desired equilibrium $\statestar$, convergence of the closed loop to a neighbourhood of $\statestar$ may still be expected, that is $\practicalset$--practical stability of the closed loop. \\
In order to show such a performance results, we assume the following:

\begin{assumption}\label{Stability and Feasibility:ass:sensitivity}
	For sets $\feasibilityset \subset \updateset \subset \statesetconstrained$ containing $\statestar$ there exist constants $\lipschitzconstantstagecost$, $\lipschitzconstantvaluefunction$, $\jumpconstantstagecost$ and $\jumpconstantvaluefunction$ such that the bounds
	\begin{align*}
		& | \stagecost(\statedisturbed, \feedbackdisturbedN(\statedisturbed, \parameterdisturbed), \parameterdisturbed) - \stagecost(\state, \feedbackN(\state, \parameter), \parameter) | \\
		& \qquad \leq \lipschitzconstantstagecost \left( \| \statedisturbed \|_{\state} + \| \parameterdisturbed \|_{\parameter}^{\max} \right) + \jumpconstantstagecost \\ \displaybreak[0]
		& | \valuefunctionN( \statedisturbed, \parameterdisturbed) - \valuefunctionN(\state, \parameter) | \\
		& \qquad \leq \lipschitzconstantvaluefunction \left( \| \statedisturbed \|_{\state} + \| \parameterdisturbed \|_{\parameter}^{\max}\right) + \jumpconstantvaluefunction
	\end{align*}
	hold for all tupels $(\statedisturbed, \state, \parameterdisturbed, \parameter)$ with $\statedisturbed \in \updateset$, $\state \in \feasibilityset$ and $\parameterdisturbed, \parameter \in \parametersetopenloop$ satisfying Assumption \ref{Setup and Preliminaries:ass:bounded disturbance}.
\end{assumption}

Based on Assumption \ref{Stability and Feasibility:ass:sensitivity} and the feasibility result given in Theorem \ref{Stability and Feasibility:thm:feasibility}, we can generalize the performance result from \cite[Theorem 3.2]{PannekGerdts2012b} to cover the case of jumps within the value function $\valuefunctionN$ due to disturbances: 

\begin{theorem}\label{Stability and Feasibility:thm:performance}
	Suppose Assumptions \ref{Setup and Preliminaries:ass:bounded disturbance} and \ref{Stability and Feasibility:ass:sensitivity} and the conditions of Theorem \ref{Stability and Feasibility:thm:feasibility} to be satisfied. Additionally let the following conditions hold:
	\begin{itemize}
		\item[(i)] A given nominal feedback $\feedbackN: \stateset \times \parameterset \to \controlset$ and a nonnegative function $\valuefunctionN: \stateset \times \parameterset \to \R_0^+$ satisfy the nominal relaxed Lyapunov inequality
		\begin{align}
			\label{Stability and Feasibility:thm:performance:eq1}
			\hspace*{-3mm}\valuefunctionN(\state, \parameter) \geq \valuefunctionN(\dynamic(\state,\feedbackN(\state, \parameter), \parameter), \parameter) + \alpha \ell(\state, \feedbackN(\state, \parameter), \parameter )
		\end{align}
		for some $\alpha \in (0, 1)$, all $\state \in \feasibilityset$ and all $\parameter \in \parameterset$.
		\item[(ii)] $\varepsilon \geq \lipschitzconstantstagecost ( \statedisturbedmax + \parameterdisturbedmax ) + \jumpconstantstagecost + ( \lipschitzconstantvaluefunction ( 2 \statedisturbedmax + 3 \parameterdisturbedmax ) + 2 \jumpconstantvaluefunction ) / \alpha$
		\item[(iii)] $\valuefunctionN(\state, \parameter) \geq \valuefunctionN(\dynamic(\state,\feedbackN(\state, \parameter), \parameter), \parameter) + \alpha \varepsilon$ holds for all $\parameter \in \parameterset$ and all $\state \in \feasibilityset \setminus \unknownset$ where $\unknownset$ denotes the minimal set which contains $\statestar$, is forward invariant with respect to all possible disturbances satisfying Assumption \ref{Setup and Preliminaries:ass:bounded disturbance}.
	\end{itemize}
	Then for the modified closed loop cost
	\begin{align*}
		\valuefunctioninftydisturbedfeedbackdisturbedN(\statedisturbed, \parameterdisturbed)   := \sum_{\closedloopindex=0}^{\infty} \overline{\stagecost}(\statedisturbed(\closedloopindex), \feedbackdisturbedN(\statedisturbed(\closedloopindex), \parameterdisturbed_{\closedloopindex}), \parameterdisturbed(\closedloopindex))
	\end{align*}
	with 
	\begin{align}
		\label{Stability and Feasibility:thm:performance:eq2}
		\overline{\stagecost}(\state, \control, \parameter) := \begin{cases}
			\max\left\{\stagecost( \state, \control, \parameter) - \varepsilon, 0\right\} & \state \in \feasibilityset \setminus \unknownset \\
			0 & \state \in \unknownset
		\end{cases}
	\end{align}
	and $\sigma :=\inf\{\valuefunctionN(\dynamic(\state, \feedbackN(\state, \parameter), \parameter), \parameter) \mid \state \in \feasibilityset \setminus \unknownset, \parameter \in \parameterset \}$ we have
	\begin{align}
		\label{Stability and Feasibility:thm:performance:eq3}
		\alpha \valuefunctioninftydisturbedfeedbackdisturbedN(\statedisturbed, \parameterdisturbed) \leq \valuefunctionN(\statedisturbed, \parameterdisturbed) - \sigma \leq \valuefunctioninfty(\statedisturbed, \parameterdisturbed) - \sigma
	\end{align}
	for all $\statedisturbed \in \updateset$.
\end{theorem}
\textbf{Proof:} Choose an arbitrary initial value $\statezero \in \feasibilityset$ and let $\closedloopindex_0 \in \N_0$ be minimal with $\state(\closedloopindex_0 + 1) \in \unknownset$ where we set $\closedloopindex_0 := \infty$ if this case does not occur.\\
	Now we reformulate \eqref{Stability and Feasibility:thm:performance:eq1} to obtain
	\begin{align*}
				\alpha \stagecost(\state, \feedbackN(\state, \parameter), \parameter) \leq \valuefunctionN(\state, \parameter) - \valuefunctionN(\dynamic(\state,\feedbackN(\state, \parameter), \parameter), \parameter).
	\end{align*}
	which allows us to incorporate the effects of disturbances and updates of the control using Assumption \ref{Stability and Feasibility:ass:sensitivity}:
	\begin{align*}
		& \alpha \stagecost(\statedisturbed(\closedloopindex), \feedbackdisturbedN(\statedisturbed(\closedloopindex), \parameterdisturbed_{\closedloopindex} ), \parameterdisturbed(\closedloopindex)) \\ \displaybreak[0]
		& \leq \alpha \stagecost(\state(\closedloopindex), \feedbackN(\state(\closedloopindex), \parameter_{\closedloopindex}), \parameter(\closedloopindex)) \\
		& \qquad + \alpha \lipschitzconstantstagecost \left( \| \statedisturbed(\closedloopindex) \|_{\state(\closedloopindex)} + \| \parameterdisturbed_{\closedloopindex} \|_{\parameter_{\closedloopindex}}^{\max} \right) + \alpha \jumpconstantstagecost \\ \displaybreak[0]
		& \leq \valuefunctionN(\state(\closedloopindex), \parameter_{\closedloopindex}) - \valuefunctionN(\dynamic(\state(\closedloopindex),\feedbackN(\state(\closedloopindex), \parameter_{\closedloopindex}), \parameter(\closedloopindex)), \parameter_{\closedloopindex}) \\
		& \qquad+  \alpha \lipschitzconstantstagecost \left( \| \statedisturbed(\closedloopindex) \|_{\state(\closedloopindex)} + \| \parameterdisturbed_{\closedloopindex} \|_{\parameter_{\closedloopindex}}^{\max}\right) + \alpha \jumpconstantstagecost \\ \displaybreak[0]
		& = \valuefunctionN(\state(\closedloopindex), \parameter_{\closedloopindex}) - \valuefunctionN(\state(\closedloopindex + 1), \parameter_{\closedloopindex}) \\
		& \qquad + \alpha \lipschitzconstantstagecost \left( \| \statedisturbed(\closedloopindex) \|_{\state(\closedloopindex)} + \| \parameterdisturbed_{\closedloopindex} \|_{\parameter_{\closedloopindex}}^{\max} \right) + \alpha \jumpconstantstagecost \\ \displaybreak[0]
		& \leq \valuefunctionN(\state(\closedloopindex), \parameter_{\closedloopindex}) - \valuefunctionN(\state(\closedloopindex + 1), \parameter_{\closedloopindex + 1}) \\
		& \qquad + \alpha \lipschitzconstantstagecost \left( \| \statedisturbed(\closedloopindex) \|_{\state(\closedloopindex)} + \| \parameterdisturbed_{\closedloopindex} \|_{\parameter_{\closedloopindex}}^{\max} \right) + \alpha \jumpconstantstagecost \\
		& \qquad + \lipschitzconstantvaluefunction \| \parameter_{\closedloopindex+1} \|_{\parameter_{\closedloopindex}}^{\max} + \jumpconstantvaluefunction
	\end{align*}
	Now we can use the feasibility result from Theorem \ref{Stability and Feasibility:thm:feasibility} to obtain
	\begin{align*}
				& \alpha \stagecost(\statedisturbed(\closedloopindex), \feedbackdisturbedN(\statedisturbed(\closedloopindex), \parameterdisturbed_{\closedloopindex} ), \parameterdisturbed(\closedloopindex)) \\ \displaybreak[0]
		& \leq \valuefunctionN(\statedisturbed(\closedloopindex), \parameterdisturbed_{\closedloopindex}) - \valuefunctionN(\statedisturbed(\closedloopindex+1), \parameterdisturbed_{\closedloopindex + 1}) \\
		& \qquad + ( \alpha \lipschitzconstantstagecost + \lipschitzconstantvaluefunction ) \left( \| \statedisturbed(\closedloopindex) \|_{\state(\closedloopindex)} + \| \parameterdisturbed_{\closedloopindex} \|_{\parameter_{\closedloopindex}}^{\max} \right) \\
		& \qquad + \lipschitzconstantvaluefunction \left( \| \statedisturbed(\closedloopindex+1) \|_{\state(\closedloopindex+1)} + \| \parameterdisturbed_{\closedloopindex+1} \|_{\parameter_{\closedloopindex+1}}^{\max} \right) \\
		& \qquad  + \lipschitzconstantvaluefunction \| \parameter_{\closedloopindex+1} \|_{\parameter_{\closedloopindex}}^{\max} + \alpha \jumpconstantstagecost + 2 \jumpconstantvaluefunction.
	\end{align*}
	Hence, using boundedness from Assumption \ref{Setup and Preliminaries:ass:bounded disturbance} reveals
	\begin{align*}
		& \alpha \stagecost(\statedisturbed(\closedloopindex), \feedbackdisturbedN(\statedisturbed(\closedloopindex), \parameterdisturbed_{\closedloopindex}), \parameterdisturbed(\closedloopindex)) \\
		& \leq \valuefunctionN(\statedisturbed(\closedloopindex), \parameterdisturbed_{\closedloopindex}) - \valuefunctionN(\statedisturbed(\closedloopindex+1), \parameterdisturbed_{\closedloopindex + 1}) \\
		& \quad + \alpha ( \lipschitzconstantstagecost ( \statedisturbedmax + \parameterdisturbedmax ) + \jumpconstantstagecost ) + \lipschitzconstantvaluefunction ( 2 \statedisturbedmax + 3 \parameterdisturbedmax ) + 2 \jumpconstantvaluefunction.
	\end{align*}
	Combining the last inequality with the condition on $\varepsilon$ in (ii), we can use condition (iii) and (i) to obtain $\valuefunctionN(\statedisturbed, \parameterdisturbed) \geq \valuefunctionN(\dynamic(\statedisturbed,\feedbackdisturbedN(\statedisturbed, \parameterdisturbed), \parameterdisturbed), \parameterdisturbed)$ for all $\state \in \feasibilityset \setminus \unknownset$ and all $(\parameter, \parameterdisturbed)$ satisfying Assumption \ref{Setup and Preliminaries:ass:bounded disturbance}. Hence, again using the bound on $\varepsilon$ and the definition of $\overline{\stagecost}$ in \eqref{Stability and Feasibility:thm:performance:eq2} we have
	\begin{align*}
		& \alpha \overline{\stagecost}(\statedisturbed(\closedloopindex), \feedbackdisturbedN(\statedisturbed(\closedloopindex), \parameterdisturbed_{\closedloopindex}), \parameterdisturbed(\closedloopindex)) \\ \displaybreak[0]
		& = \max\left\{ \alpha \stagecost(\statedisturbed(\closedloopindex), \feedbackdisturbedN(\statedisturbed(\closedloopindex), \parameterdisturbed_{\closedloopindex}), \parameterdisturbed(\closedloopindex)) - \alpha \varepsilon, 0\right\} \\ \displaybreak[0]
		& \leq \valuefunctionN(\statedisturbed(\closedloopindex), \parameterdisturbed_{\closedloopindex}) - \valuefunctionN(\statedisturbed(\closedloopindex+1), \parameterdisturbed_{\closedloopindex+1}).
	\end{align*}
	For $\closedloopindex \geq \closedloopindex_0 + 1$ the invariance of $\unknownset$ gives us $\statedisturbed(\closedloopindex) \in \unknownset$ and hence $\overline{\stagecost}(\statedisturbed(\closedloopindex), \feedbackdisturbedN(\statedisturbed(\closedloopindex), \parameterdisturbed_{\closedloopindex}), \parameterdisturbed(\closedloopindex)) = 0$. Additionally, since $\sigma$ is the minimal cost after entry in $\unknownset$, we have $\valuefunctionN(\statedisturbed(\closedloopindex), \parameterdisturbed_{\closedloopindex}) \geq \sigma$ for all $\closedloopindex \leq \closedloopindex_0 + 1$. Now the feasibility result from Theorem \ref{Stability and Feasibility:thm:feasibility} allows us to sum the stage costs over $\closedloopindex$ and obtain
	\begin{align*}
		& \alpha \sum_{\closedloopindex = 0}^{K} \overline{\stagecost}(\statedisturbed(\closedloopindex), \feedbackdisturbedN(\statedisturbed(\closedloopindex), \parameterdisturbed_{\closedloopindex}), \parameterdisturbed(\closedloopindex)) \\ \displaybreak[0]
		& = \alpha \sum_{\closedloopindex = 0}^{K_0} \overline{\stagecost}(\statedisturbed(\closedloopindex), \feedbackdisturbedN(\statedisturbed(\closedloopindex), \parameterdisturbed_{\closedloopindex}), \parameterdisturbed(\closedloopindex)) \\ \displaybreak[0]
		& \leq \valuefunctionN(\statedisturbed(0), \parameterdisturbed_{0}) - \valuefunctionN(\statedisturbed(K_0 + 1), \parameterdisturbed_{K_0 + 1}) \\ \displaybreak[0]
		& \leq \valuefunctionN(\statedisturbed(0), \parameterdisturbed_{0}) - \sigma.
	\end{align*}
	where $K_0 := \min \{ K, \closedloopindex_0 \}$. Using that $K \in \N$ was arbitrary we can conclude that $(\valuefunctionN(\statedisturbed(0), \parameterdisturbed_{0}) - \sigma)/\alpha$ is an upper bound for $\valuefunctioninftydisturbedfeedbackdisturbedN(\statedisturbed(0), \parameterdisturbed)$. Last, using arbitraryness of the initial value $\statezero \in \feasibilityset$, assertion \eqref{Stability and Feasibility:thm:performance:eq3} follows. \qed
\begin{remark}
	Within Theorem \ref{Stability and Feasibility:thm:performance} $\unknownset$ is implicitly defined. For approximation techniques of this set we refer to \cite{GrimmMessinaTeel2005, GrueneRantzer2008}.
\end{remark}
Last, we can use the previous performance result to show $\mathcal{P}$--practical asymptotic stability of system \eqref{Setup and Preliminaries:eq:plant}.
\begin{theorem}\label{Stability and Feasibility:thm:stability}
	Suppose the conditions of Theorem \ref{Stability and Feasibility:thm:performance} hold and additionally there exist $\cK_\infty$ functions $\alpha_1$, $\alpha_2$ such that
	\begin{align}
		\label{Stability and Feasibility:thm:stability:eq1}
		\alpha_1( \| \statedisturbed \|_{\statestar} + \| \parameterdisturbed \|_{\parameterstar}^{\max} ) & \leq \valuefunction(\statedisturbed, \parameter) \leq  \alpha_2( \| \statedisturbed \|_{\statestar} + \| \parameterdisturbed \|_{\parameterstar}^{\max} )
	\end{align}
	holds for all $\statedisturbed \in \updateset \setminus \mathcal{P}$ with $\mathcal{P} = \unknownset$. Then $\statestar$ is $\mathcal{P}$--practically asymptotically stable on $\updateset$.
\end{theorem}
\textbf{Proof:} The definition of $\overline{\stagecost}$ in \eqref{Stability and Feasibility:thm:performance:eq2} and the property $\valuefunction(\statedisturbed, \parameterdisturbed) \geq \valuefunctionN(\dynamic(\statedisturbed,\feedbackdisturbedN(\statedisturbed, \parameterdisturbed), \parameterdisturbed), \parameterdisturbed)$ shown in the proof of Theorem \ref{Stability and Feasibility:thm:performance} guarantee the conditions of Theorem \ref{Setup and Preliminaries:thm:p practical} to hold showing the assertion. \qed

\begin{remark}
	(i) To cover the case of $\mathcal{P}$--practically stable nominal systems \eqref{Setup and Preliminaries:eq:control system} condition \eqref{Stability and Feasibility:thm:performance:eq1} can be relaxed to
	\begin{align*}
		& \valuefunctionN(\state, \parameter) - \valuefunctionN(\dynamic(\state,\feedbackN(\state, \parameter), \parameter), \parameter) \\
		& \geq \min \{ \alpha \ell(\state, \feedbackN(\state, \parameter), \parameter ) - \overline{\varepsilon}, \ell(\state, \feedbackN(\state, \parameter), \parameter ) - \overline{\varepsilon} \}
	\end{align*}
	as shown in \cite{GrueneRantzer2008}. Then the performance and stability results from Theorems \ref{Stability and Feasibility:thm:performance} and \ref{Stability and Feasibility:thm:stability} hold if we consider the modified lower bound $\varepsilon \geq \lipschitzconstantstagecost ( \statedisturbedmax + \parameterdisturbedmax ) + \jumpconstantstagecost + ( \lipschitzconstantvaluefunction ( 2 \statedisturbedmax + 3 \parameterdisturbedmax ) + 2 \jumpconstantvaluefunction + \overline{\varepsilon} ) / \alpha$. \\
	(ii) Using identical arguments as in the proof of Theorem \ref{Stability and Feasibility:thm:stability}, $m$--step MPC control laws can be handled. To this end, condition \eqref{Stability and Feasibility:thm:performance:eq1} can replaced by
	\begin{align*}
		& \valuefunctionN(\state, \parameter) - \valuefunctionN(\state_{\control, \parameter}(m; \state), \parameter) \\
		& \geq \alpha \sum_{\openloopindex=0}^{m-1} \ell(\state_{\control, \parameter}(\openloopindex; \state), \feedbackN(\state_{\control, \parameter}(\openloopindex; \state), \parameter), \parameter )
	\end{align*}
	and both $\overline{\stagecost}(\state, \control, \parameter)$ and $\sigma$ need to be redefined accordingly.
\end{remark}

Here, we like to stress that the stability result of Theorem \ref{Stability and Feasibility:thm:stability} holds for any updated feedback law $\feedbackdisturbedN$ satisfying Assumption \ref{Stability and Feasibility:ass:sensitivity}. In the next section, we verify these requirements for particular updating strategies known from the MPC literature.

\section{Update Techniques}
\label{Section:Update Techniques}

Currently, there are two main lines of updating techniques used in MPC which are heavily influenced by nonlinear optimization theory: For one, as shown in \cite{ZavalaBiegler2009}, sensitivity information \cite{Fiacco1983} can be applied in a similar manner as in offline optimal control, cf., e.g. \cite{GroetschelKrumkeRambau2001}, to update the control law based on newly arrived information. Different from that, interim updates on optimality or feasibility can be employed to cope with measureable disturbances. The key idea here is that the solutions of two consecutive MPC iterates are close to each other. Hence, for the realtime iteration approach \cite{DiehlBockSchloeder2005}, one can show that typically one iteration of the underlining nonlinear optimization routine is sufficient to guarantee stability. Extending this idea, hierachical methods \cite{BockDiehlKostinaSchloeder2007} may operate on different levels and different time scales of the problem to robustify the MPC method.\\
While the ideas of these updating methods are simple, corresponding proofs are quite involved and require a lot of notation used in nonlinear optimization theory. Since we make use of respective properties of these methods only, we condense the additional notation to a minimum and refer to the original publications \cite{BockDiehlKostinaSchloeder2007, DiehlBockSchloeder2005, Fiacco1983} for details.

\subsection{Sensitivity Analysis}

Given the advanced step setting of Algorithm \ref{Setup and Preliminaries:alg:mpc}, it is possible to solve a nominal optimal control problem ahead of time and to additionally compute sensitivity information $\frac{\partial \controlstar}{\partial \state}$, $\frac{\partial \controlstar}{\partial \parameter}$. Hence, upon arrival of new measurements of $\statedisturbed$ and $\parameterdisturbed$, the optimal control can be updated via
\begin{align}
	\label{Sensitivity:eq:update formula}
	\hspace*{-2mm}\controldisturbed(\cdot, \statedisturbed, \parameterdisturbed) := \controlstar(\cdot, \state, \parameter) + %
	\begin{pmatrix}
		\frac{\partial \controlstar}{\partial \state}(\cdot, \state, \parameter) \\
		\frac{\partial \controlstar}{\partial \parameter}(\cdot, \state, \parameter)
	\end{pmatrix}^\top%
	\begin{pmatrix}
		\statedisturbed - \state \\
		\parameterdisturbed - \parameter
	\end{pmatrix}
\end{align}
which allows us to set $\feedbackdisturbedN(\statedisturbed, \parameterdisturbed) = \controldisturbed(0, \statedisturbed, \parameterdisturbed)$. A proof under which conditions such an update reveals an optimal control is given in the following theorem from \cite{Fiacco1983} where we adapted the notation to match the considered MPC case:
\begin{theorem}\label{Sensitivity:thm:sensitivity}
	Consider $\dynamic$ and $\stagecost$ to be twice continuously differentiable in a neighbourhood $\updateset \subset \statesetconstrained$ of the nominal solution $\controlstar(\cdot, \state, \parameter)$ of Problem \eqref{Setup and Preliminaries:eq:value function}. If the linear independence constraint qualification (LICQ), the second order sufficient optimality conditions (SSOC) and the strict complementarity condition (SCC) hold in this neighbourhood, then
	\begin{itemize}
		\item $\controlstar(\cdot, \state, \parameter)$ is an isolated local minimizer and the respective Lagrange multipliers are unique,
		\item there exists a unique local minimizer $\controlstar(\cdot, \statedisturbed, \parameterdisturbed)$ for $(\statedisturbed, \parameterdisturbed )$ in a neighbourhood of $( \state, \parameter )$ which satisfies LICQ, SSOC and SCC and is differentiable with respect to $\statedisturbed$ and $\parameterdisturbed$,
		\item there exist a Lipschitz constant $\lipschitzconstantvaluefunction$ such that
		\begin{align}
			\label{Sensitivity:thm:sensitivity:eq1}
			& | \valuefunctionN( \statedisturbed, \parameterdisturbed) - \valuefunctionN(\state, \parameter) | \leq \lipschitzconstantvaluefunction \left( \| \statedisturbed \|_{\state} + \| \parameterdisturbed \|_{\parameter}^{\max} \right)
		\end{align}
		holds and
		\item there exist a Lipschitz constant $L_\control$ such that for the updated control $\controldisturbed(\cdot, \statedisturbed, \parameterdisturbed)$ from \eqref{Sensitivity:eq:update formula} the following estimate holds:
		\begin{align*}
			\| \controldisturbed(\cdot, \statedisturbed, \parameterdisturbed) \|_{\controlstar(\cdot, \state, \parameter)}^{\max} \leq \lipschitzconstantcontrol \left( \| \statedisturbed \|_{\state} + \| \parameterdisturbed \|_{\parameter}^{\max} \right)
		\end{align*}
	\end{itemize}
\end{theorem}

Now, the following result allows us to conclude that Assumption \ref{Stability and Feasibility:ass:sensitivity} holds if the conditions of Assumption \ref{Setup and Preliminaries:ass:bounded disturbance} and Theorem \ref{Sensitivity:thm:sensitivity} apply.
\begin{proposition}\label{Sensitivity:prop:assumptions}
	If there exist sets $\feasibilityset \subset \updateset \subset \statesetconstrained$ containing $\statestar$ such that $\feasibilityset$ is control forward invariant and the conditions of Theorem \ref{Sensitivity:thm:sensitivity} hold for all tupels $(\statedisturbed, \state, \parameterdisturbed, \parameter)$ with $\statedisturbed \in \updateset$, $\state \in \feasibilityset$ and $\parameterdisturbed, \parameter \in \parametersetopenloop$ satisfying Assumption \ref{Setup and Preliminaries:ass:bounded disturbance}, then Assumption \ref{Stability and Feasibility:ass:sensitivity} holds.
\end{proposition}
\textbf{Proof:} The conclusion follows directly from \eqref{Sensitivity:thm:sensitivity:eq1} and the fact that differentiability of $\stagecost$ implies the existence of a local Lipschitz constant $\lipschitzconstantstagecost$. \qed

Hence, we can conclude the following:

\begin{conclusion}\label{Sensitivity:con:stability}
	Suppose Assumption \ref{Setup and Preliminaries:ass:bounded disturbance} and the conditions of Theorem \ref{Stability and Feasibility:thm:feasibility} and Proposition \ref{Sensitivity:prop:assumptions} to hold. If there exist $\cK_\infty$ functions $\alpha_1$, $\alpha_2$ such that \eqref{Stability and Feasibility:thm:stability:eq1} is true for all $\statedisturbed \in \updateset \setminus \mathcal{P}$ with $\mathcal{P} = \unknownset$, then the control computed by Algorithm \ref{Setup and Preliminaries:alg:mpc} $\mathcal{P}$--practically asymptotically stabilizes $\statestar$ on $\updateset$.
\end{conclusion}
\textbf{Proof:} Assumption \ref{Setup and Preliminaries:ass:bounded disturbance} and Proposition \ref{Sensitivity:prop:assumptions} guarantee the conditions of Theorem \ref{Stability and Feasibility:thm:stability} to hold showing the assertion. \qed

\subsection{Realtime Iterations and Hierarchical MPC}

In contrast to sensitivity based updates, realtime iterations can be performed without precomputing sensitivity matrices. Instead, only a single Newton step of the optimization routine is executed and the updated control is applied right away. It has been shown in \cite[Theorem 4.1]{DiehlBockSchloeder2005}, that such a procedure leads to a contraction of the disturbance around the nominal solution and the updated solution converges towards the stationary point of the disturbed problem. Considering the notation introduced in this paper, this results reads as follows:
\begin{theorem}\label{Realtime:thm:contraction}
	Consider $\dynamic$ and $\stagecost$ to be twice continuously differentiable and the approximation of the second order derivative of the Lagrangian to be continuous with bounded inverse in a neighbourhood $\updateset \subset \statesetconstrained$ of the nominal solution $\controlstar(\cdot, \state, \parameter)$ of Problem \eqref{Setup and Preliminaries:eq:value function}. Additionally suppose that both the approximation error of the Lagrangian for a linear interpolation of the disturbance and the disturbance impact on the approximation of the Lagrangian are bounded. If the disturbance is sufficiently small and the disturbed state trajectory $\statedisturbed_{\control, \parameterdisturbed}(\cdot; \statedisturbed)$ is feasible, then applying one Newton step we obtain
	\begin{itemize}
		\item[(i)] feasibility of the updated state trajectory $\statedisturbed_{\controldisturbed, \parameterdisturbed}(\cdot; \statedisturbed)$ and
		\item[(ii)] existence of constants $\realtimeerrorconstant < 1$ and $\realtimeerrorlinear < \infty$ such that the following error bound is guaranteed:
		\begin{align}
			\label{Realtime:thm:contraction:eq1}
			\| \statedisturbed_{\controldisturbed, \parameterdisturbed}(\cdot; \statedisturbed) \|_{\state_{\control, \parameter}(\cdot; \state)}^{\max} \leq & \; \realtimeerrorconstant \| \statedisturbed_{\control, \parameterdisturbed}(\cdot; \statedisturbed) \|_{\state_{\control, \parameter}(\cdot; \state)}^{\max} \\
			& + \realtimeerrorlinear \left( {\| \statedisturbed_{\control, \parameterdisturbed}(\cdot; \statedisturbed) \|_{\state_{\control, \parameter}(\cdot; \state)}^{\max}} \right)^2.\nonumber 
		\end{align}
	\end{itemize}
	Additionally, the Newton sequence converges to the exact stationary point of the disturbed problem.
\end{theorem}

Based on this result and Assumption \ref{Setup and Preliminaries:ass:bounded disturbance}, we can show that our fundamental Assumption \ref{Stability and Feasibility:ass:sensitivity} is satisfied:

\begin{proposition}\label{Realtime:prop:assumptions}
	If there exist sets $\feasibilityset \subset \updateset \subset \statesetconstrained$ containing $\statestar$ such that $\feasibilityset$ is control forward invariant and the conditions of Theorem \ref{Realtime:thm:contraction} hold for all tupels $(\statedisturbed, \state, \parameterdisturbed, \parameter)$ with $\statedisturbed \in \updateset$, $\state \in \feasibilityset$ and $\parameterdisturbed, \parameter \in \parametersetopenloop$ satisfying Assumption \ref{Setup and Preliminaries:ass:bounded disturbance}, then Assumption \ref{Stability and Feasibility:ass:sensitivity} holds.
\end{proposition}
\textbf{Proof:} Since $\dynamic$ is Lipschitz, boundedness of the state and parameter disturbances given by Assumption \ref{Setup and Preliminaries:ass:bounded disturbance} allows us to conclude that there exists a constant $\lipschitzconstantstate$ such that the quadratic bound \eqref{Realtime:thm:contraction:eq1} can be overbounded linearly by
\begin{align*}
	\| \statedisturbed_{\controldisturbed, \parameterdisturbed}(\cdot; \statedisturbed) \|_{\state_{\control, \parameter}(\cdot; \state)}^{\max} \leq \lipschitzconstantstate \left( \| \statedisturbed \|_{\state} + \| \parameterdisturbed \|_{\parameter}^{\max} \right).
\end{align*}
Moreover, since the Newton iterates are contracting and therefore bounded, we obtain that there exists a constant $\lipschitzconstantcontrol$ such that
\begin{align*}
	\| \controldisturbed(\cdot, \statedisturbed, \parameterdisturbed) \|_{\controlstar(\cdot, \state, \parameter)}^{\max} \leq \lipschitzconstantcontrol \left( \| \statedisturbed \|_{\state} + \| \parameterdisturbed \|_{\parameter}^{\max} \right)
\end{align*}
holds. Now, due to continuity of $\dynamic$ and $\stagecost$ boundedness of the disturbed state trajectory $\statedisturbed_{\controldisturbed, \parameterdisturbed}(\cdot; \statedisturbed)$ and the updated control $\controldisturbed(\cdot, \statedisturbed, \parameterdisturbed)$ allows us to conclude that there exists a constant $\lipschitzconstantvaluefunction$ such that \eqref{Sensitivity:thm:sensitivity:eq1} holds. Similar to the proof of Proposition \ref{Sensitivity:prop:assumptions}, differentiability of $\stagecost$ implies the existence of a local Lipschitz constant $\lipschitzconstantstagecost$ showing the assertion. \qed

Hence, $\mathcal{P}$--practical asymptotic stability of the closed loop can be concluded similar to the sensitivity based scenario:

\begin{conclusion}\label{Realtime:con:stability}
	Consider Assumption \ref{Setup and Preliminaries:ass:bounded disturbance} and the conditions of Theorem \ref{Stability and Feasibility:thm:feasibility} and Proposition \ref{Realtime:prop:assumptions} to apply. If there exist $\cK_\infty$ functions $\alpha_1$, $\alpha_2$ such that \eqref{Stability and Feasibility:thm:stability:eq1} holds for all $\statedisturbed \in \updateset \setminus \mathcal{P}$ with $\mathcal{P} = \unknownset$, then the control computed by Algorithm \ref{Setup and Preliminaries:alg:mpc} $\mathcal{P}$--practically asymptotically stabilizes $\statestar$ on $\updateset$.
\end{conclusion}
\textbf{Proof:} Similar to the proof of Conclusion \ref{Sensitivity:con:stability}. \qed

Here, we like to stress that the convergence result for hierarchically structured MPC updates show the properties (i) and (ii) of Theorem \ref{Realtime:thm:contraction}, cf. \cite[Theorem 6]{BockDiehlKostinaSchloeder2007}. Hence, the results of Proposition \ref{Realtime:prop:assumptions} and Conclusion \ref{Realtime:con:stability} also apply in the context of hierarchical MPC updates.

\begin{remark}
	(i) Note that the feasibility and update sets $\feasibilityset$ and $\updateset$ within Theorems \ref{Stability and Feasibility:thm:performance}, \ref{Sensitivity:thm:sensitivity}, \ref{Realtime:thm:contraction} and Propositions \ref{Sensitivity:prop:assumptions}, \ref{Realtime:prop:assumptions} are not necessarily large. This become clear by the requirement of identical open loop control structures in Theorem \ref{Sensitivity:thm:sensitivity} which is typically only satisfied in a small neighbourhood. Hence, both the maximal disturbances $\statedisturbedmax$, $\parameterdisturbedmax$ and the Lipschitz constants $\lipschitzconstantvaluefunction$, $\lipschitzconstantstagecost$ are comparably small. Since the properties are required only locally, i.e. for one MPC step, the closed loop control structure may still vary. The closed loop feasibility and update sets and constants can then be obtained by the union of the local sets and maximum of local values respectively. \\
	(ii) Within the proofs in the section the jump condition constants $\jumpconstantvaluefunction$ and $\jumpconstantstagecost$ in Assumption \ref{Stability and Feasibility:ass:sensitivity} can be set to zero. Consequently, the results from Section \ref{Section:Stability and Feasibility} apply to a wider range of updates than the ones we considered here.
\end{remark}

\section{Numerical Example}
\label{Section:Numerical Results}

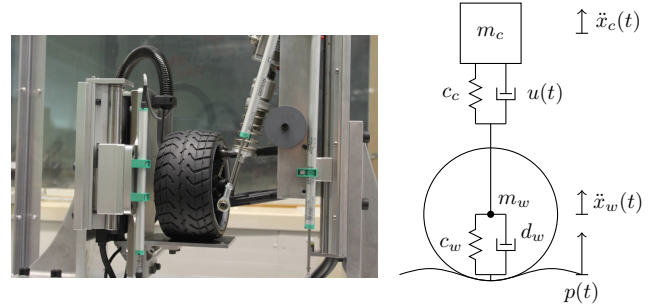
\begin{figure}[!ht]
	\begin{minipage}{0.27\textwidth}
		\includegraphics[width=0.98\textwidth]{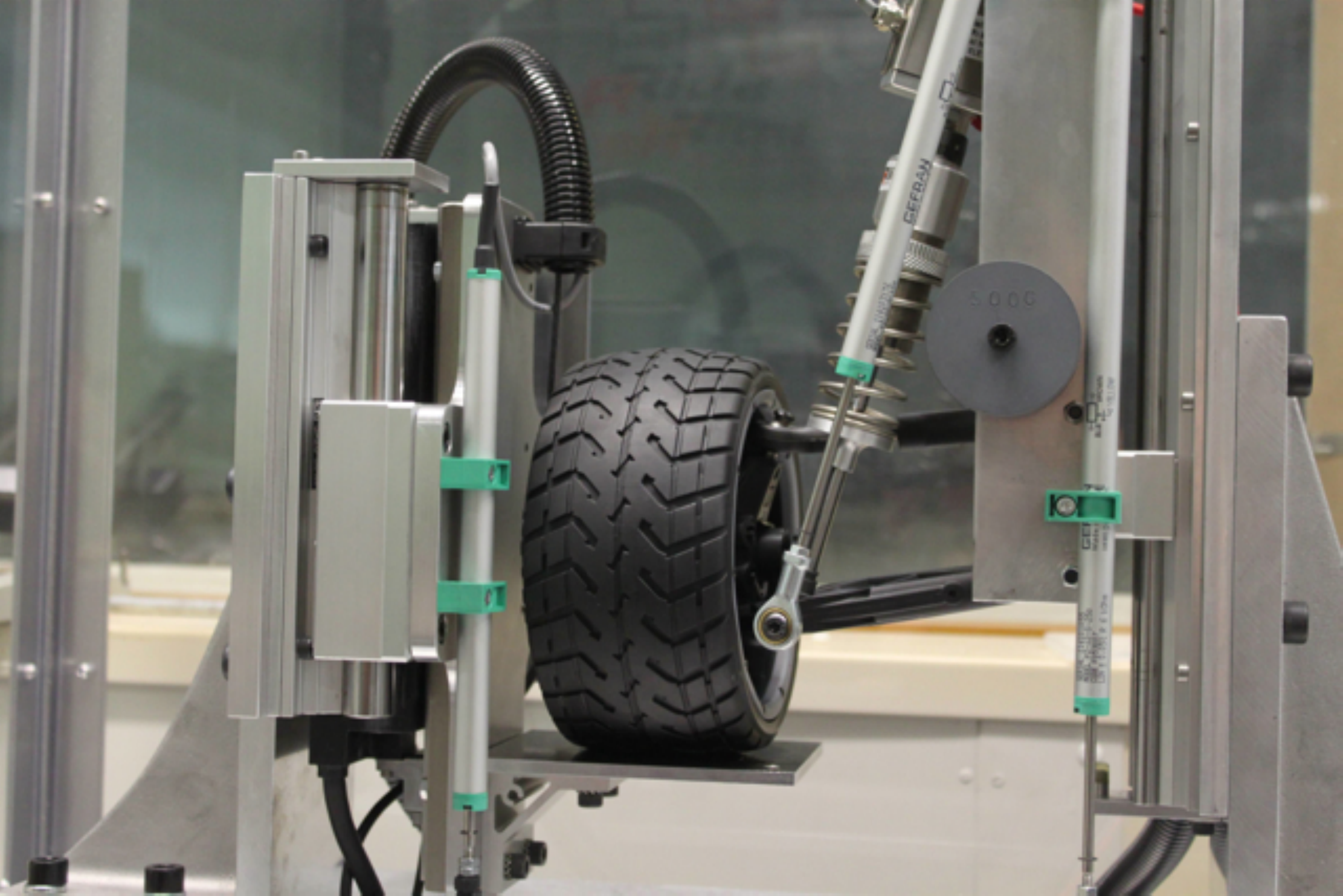}	
	\end{minipage}
	\begin{minipage}{0.19\textwidth}
	\begin{center}
    	\begin{tikzpicture}[scale = 0.8, every node/.style={scale=0.8}]
        \draw (4,4) rectangle (5,5);
        \draw (4.5,4.5) node {$m_c$};
       
        \draw (4.75,3) -- (4.75,3.4);
        \draw (4.75,3.4) -- (4.6,3.4) -- (4.6,3.6);
        \draw (4.75,3.4) -- (4.9,3.4) -- (4.9,3.6);
        \draw (4.65,3.5) -- (4.85,3.5);
        \draw (4.75,3.5) -- (4.75,4);
        \draw (5.4,3.5) node {$u(t)$};
        
        \draw (4.25,3) -- (4.25,3.25) -- (4.15,3.3) -- (4.35,3.4) -- (4.15,3.5) -- (4.35,3.6) -- (4.15,3.7) -- (4.25,3.75) -- (4.25,4);
        \draw (3.85,3.5) node {$c_c$};
        
        \draw (4.25,3) -- (4.75,3);
        \draw (4.5,3) -- (4.5,1.5);
        \draw (4.5,1.5) circle (1.1);
        \draw[fill] (4.5,1.5) circle (0.05) node[above right] {$m_w$};
        \draw (4.25,1.5) -- (4.75,1.5);
        
        \draw (4.75,0.5) -- (4.75,0.9);
        \draw (4.75,0.9) -- (4.6,0.9) -- (4.6,1.1);
        \draw (4.75,0.9) -- (4.9,0.9) -- (4.9,1.1);
        \draw (4.65,1.0) -- (4.85,1.0);
        \draw (4.75,1.0) -- (4.75,1.5);
        \draw (5.2,1.2) node {$d_w$};
        
        \draw (4.25,0.5) -- (4.25,0.75) -- (4.15,0.8) -- (4.35,0.9) -- (4.15,1) -- (4.35,1.1) -- (4.15,1.2) -- (4.25,1.25) -- (4.25,1.5);
        \draw (3.85,1) node {$c_w$};
                
        \draw (4.25,0.5) -- (4.75,0.5);
        \draw (4.5,0.5) -- (4.5,0.4);
        \draw (3,0.5) sin (3.5,0.6) cos (4,0.5) sin (4.5,0.4) cos (5,0.5) sin (5.5,0.6) cos (6,0.5);
        \draw (6,0.2) node {$\parameter(t)$};
        
         \draw[->] (6,1.5) -- (6,1.9);
         \draw (5.9,1.5) -- (6.1,1.5);
         \draw (6.6,1.7) node {$\ddot{\state}_w(t)$};

         \draw[->] (6,4.5) -- (6,4.9);
         \draw (5.9,4.5) -- (6.1,4.5);
         \draw (6.6,4.7) node {$\ddot{\state}_c(t)$};
		\draw[->] (6,0.5) -- (6,1.2);
		\draw (5.9,0.5) -- (6.1,0.5);
      	\end{tikzpicture}
      	\end{center}
  	\end{minipage}
    \caption{Test bench and schematic drawing of the quartercar}
    \label{Numerical Results:fig:quartercar model}
\end{figure}

The previous analysis was motivated by considering a quarter car application with forward sensors for measuring road data and state sensors for the position of both wheel and chassis, cf. Figure \ref{Numerical Results:fig:quartercar model}.

The dynamics of the quarter car model are given by the set of second order ordinary differential equations
\begin{align} 
 m_w \ddot{\state}_w(t) &= c_c(\state_c(t) \hspace*{-0.5mm} - \hspace*{-0.5mm} \state_w(t)) \hspace*{-0.5mm} + \hspace*{-0.5mm} \control(t)(\dot{\state}_c(t)\hspace*{-0.5mm}-\hspace*{-0.5mm}\dot{\state}_w(t)) \nonumber\\
	    & \quad - \hspace*{-0.5mm} c_w(\state_w(t) \hspace*{-0.5mm} - \hspace*{-0.5mm} \parameter(t)) \hspace*{-0.5mm} - \hspace*{-0.5mm} d_w(\dot{\state}_w(t) \hspace*{-0.5mm} - \hspace*{-0.5mm} \dot{\parameter}(t))\label{NumericalResults:model}\\  
 m_c \ddot{\state}_c(t) &= -\hspace*{-0.5mm} c_c(\state_c(t) \hspace*{-0.5mm} - \hspace*{-0.5mm} \state_w(t)) \hspace*{-0.5mm} - \hspace*{-0.5mm} \control(t)(\dot{\state}_c(t) \hspace*{-0.5mm} - \hspace*{-0.5mm} \dot{\state}_w(t)) \nonumber
\end{align}
where $\state_w$ and $\state_c$ denote the state of the wheel and the chassis and $m_w$ and $m_c$ the respective wheel and chassis masses. The parameters $d_w$ and $c_w$ represent the spring and the damper coefficients modelling the wheel to ground interaction whereas $c_c$ denotes the spring constant of the chassis, see Table \ref{Numerical Results:tab:quartercar} for respective values.
\begin{table}[!htb]
	\begin{center}
		\begin{tabular}{|c|c|c|c|} \hline
			Name & Symbol & Quantity & Unit\\ \hline\hline
			mass wheel & $m_w$ & $35$ & $kg$ \\
			mass chassis & $m_c$ & $325$ & $kg$\\
			spring constant wheel & $c_w$ & $0.2$ & $kN/m$\\
			damper constant wheel & $d_w$ & $150$ & $kNs/m$\\
			spring constant chassis & $c_c$ & $20$ & $kN/m$\\
		\hline
		\end{tabular}
		\caption{Parameters for the quarter car example from \cite{RettigStryk2005}}
		\label{Numerical Results:tab:quartercar}
	\end{center}
\end{table}
Moreover, the functions $\control(t)$ and $\parameter(t)$ represent the controllable damper constant of the chassis and the road profile. Note that within system \eqref{NumericalResults:model}, the states $\state_w$ and $\state_c$ do not correspond to the actual heights above ground of the center of masses of both wheel and body. Instead, the equilibrium of the system has been shifted to the origin. In the following, we use our standard notation and abbreviate the system state via $\state(t) := \left[\state_w, \state_c\right]^{\top}$. 

\begin{center}
	\begin{figure}[!htb]
		\includegraphics[width=0.47\textwidth]{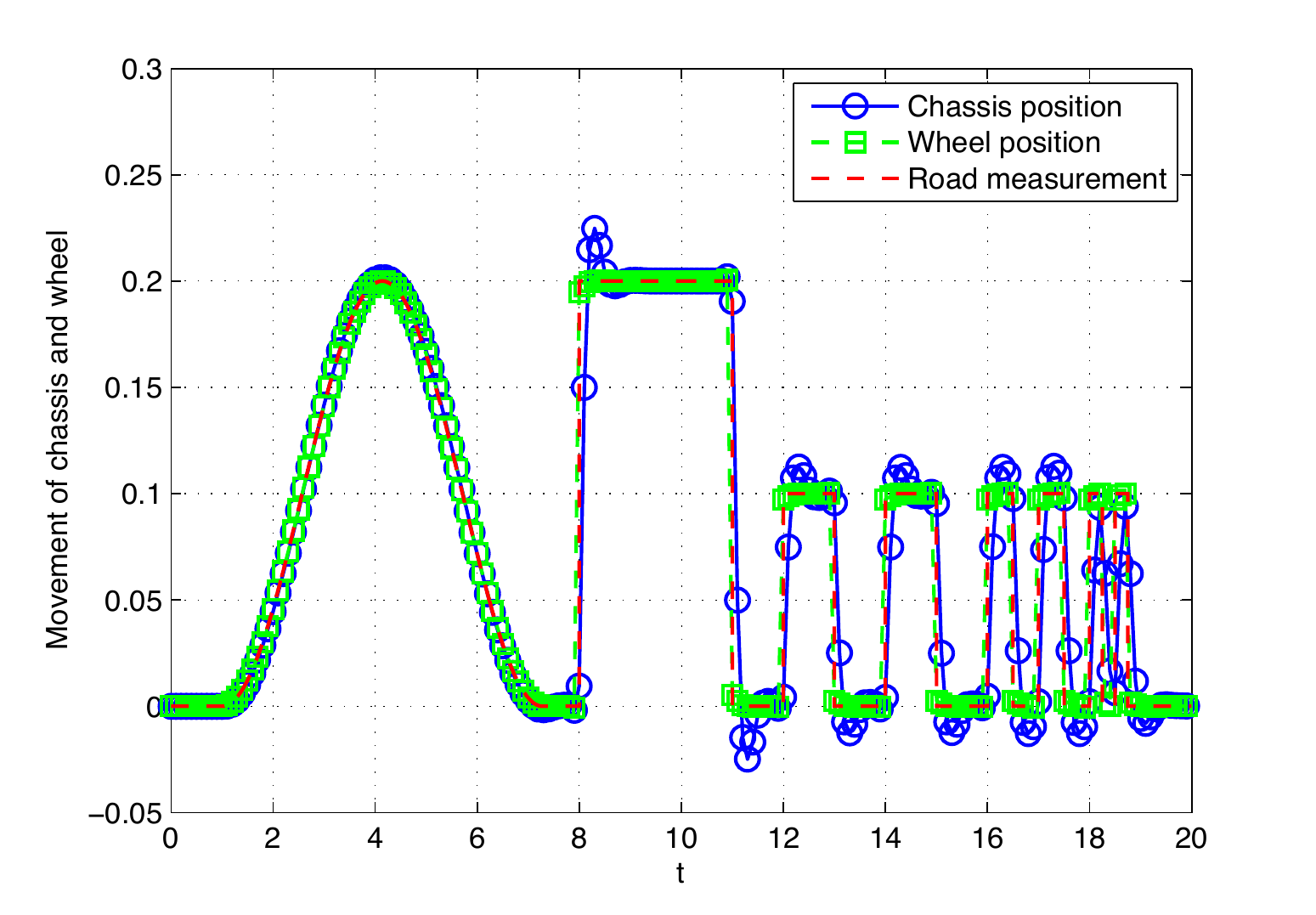}
		\caption{Movement of chassis and wheel in $m$ using the nominal MPC law without disturbances}
		\label{Numerical Results:fig:undisturbed movement}
	\end{figure}
\end{center}

Within the MPC setup we considered the criteria comfort and safety in accordance with ISO/IEC 2631 \cite{ISO2631}. In particular, the jerk of the chassis mass was employed in
\begin{equation*}
 \costfunction_{\text{comfort}}(\state(t),\control(t),\parameter(t)) = \int_0^{T} \left( \dddot{\state}_c \right)^2dt,
\end{equation*}
to measure the comfort of the driver. To measure safety, the functional
\begin{align*}
 \costfunction_{\text{safety}}(\state(t),\control(t),\parameter(t)) = \int_0^{T} ( & c_w ( \state_w(t) - \parameter(t) ) \\
  & + d_w (\dot{\state}_w(t) - \dot{\parameter}(t) )^2 dt.
\end{align*}
is used which represent the deviation from the equilibrium tire--surface force. 
In order to adapt to different driving situations and driver adjustments, these criteria are weighted by parameters $\mu_{\text{comfort}}$ and $\mu_{\text{safety}}$
\begin{align}
	& \costfunctionN(\state,\control,\parameter) = \\
 	& \quad = \sum_{\openloopindex=0}^{\horizon-1} \mu_{\text{comfort}} \int_{\openloopindex T}^{(\openloopindex+1)T} \costfunction_{\text{comfort}}(\state(t),\control(t),\parameter(t)) dt \nonumber \\
 	& \qquad + \mu_{\text{safety}} \int_{\openloopindex T}^{(\openloopindex+1)T} \costfunction_{\text{safety}}(\state(t),\control(t),\parameter(t)) dt. \nonumber
\end{align}
For our computations, we considered the control value set $\controlsetconstrained = \left[0.5, 3.0\right]$ and a sampling period $T = T_\control = 0.1s$ during which the control is held constant. Measurements of the road data, however, are taken at a sampling period of $T_\parameter = 0.002s$ from an artificial test track, cf. the dashed line in Figure \ref{Numerical Results:fig:undisturbed movement}. Note that since we utilized a point model of the tire, we focused on large road excitations assuming that a more realistic tire model such as CDTire \cite{BGHT2011} reduces the effect of small excitations.\\
In order to solve the resulting optimal control problem, we set $\horizon = 5$, $\mu_{\text{comfort}} = 10$, $\mu_{\text{safety}} = 1$ and applied a direct approach. To this end, we discretized the problem using the sampling period of the control $T_\control$ and set the optimization tolerance of the used SQP method to $10^{-6}$. Within each sampling period, the dynamics and the cost functional are evaluated via the DoPri5 method with error tolerance $10^{-6}$. Last, the required twice continuously differentiable road profile $\parameter(t)$ is obtained by means of a Fast Fourier Transformation (FFT) for the nominal measurements over the optimization horizon. The FFT itself is based on the $\horizon \cdot T_\control / T_\parameter + 1 = 251$ road measurements contained within the open loop optimization horizon and recomputed for each sampling point of the control. Here, we like to mention that although the coefficients need to be calculated only once for each MPC step, the evaluation of the FFT interpolation polynome remains the computationally most expensive part since it is required by the differential equation solver. To reduce this additional effort, the integration method can be restricted to operate on the sampling instances of the road measurements only. Alternatively, since a more realistic tire model compensates small excitations with high frequencies, a low pass filter for the FFT may be used to reduce the number of coefficients.

\begin{center}
	\begin{figure}[!htb]
		\includegraphics[width=0.47\textwidth]{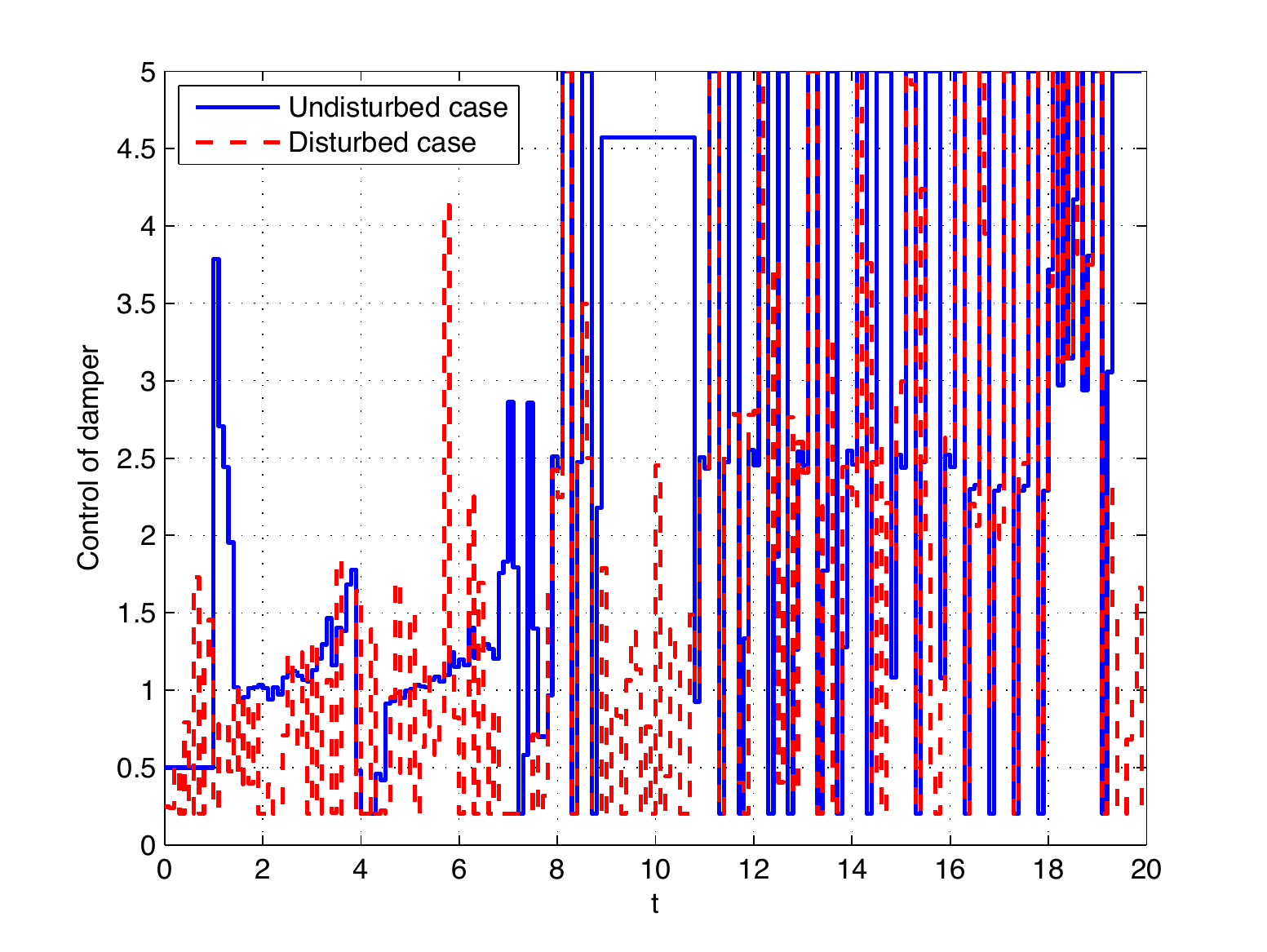}
		\caption{Computed damping strategy for given road excitation in $kNs/m$}
		\label{Numerical Results:fig:undisturbed and disturbed control}
	\end{figure}
\end{center}

As one observes from Figure \ref{Numerical Results:fig:undisturbed movement}, the $\state_w$ and $\state_c$ trajectories follow the road excitation profile nicely. Additionally, the expected overshoot of the chassis fades out within a few sampling instances due to the computed control strategy shown in Figure \ref{Numerical Results:fig:undisturbed and disturbed control}.\\
\vspace{-0.5cm}
\vspace{-0.5cm}
To incorporate disturbances, we first introduced measurement errors of the road profile due to sensor scattering. Here, we considered the error to be uniformly distributed in the interval $[-0.005, 0.005]$ which corresponds to the typical accuracy of a laser scanner of $5 mm$ in our simulated application. Regarding Assumption \ref{Setup and Preliminaries:ass:bounded disturbance} we therfore obtain the parameter error bound $\parameterdisturbedmax = 0.005$. While the closed loop trajectories for both chassis and wheel are almost identical to the undisturbed case, the control strategy differs significantly from the undisturbed one, cf. Figure \ref{Numerical Results:fig:undisturbed and disturbed control}.
\vspace{-0.3cm}
\vspace{-0.5cm}
In this example, it is noteworthy that for the sinoidal excitation the boundary conditions on the control are now active for many MPC iterates. Hence, conditions of Theorem \ref{Sensitivity:thm:sensitivity} do not hold true and therefore applicability of the sensitivity based update approach is questionable.

\begin{center}
	\begin{figure}[!htb]
		\includegraphics[width=0.47\textwidth]{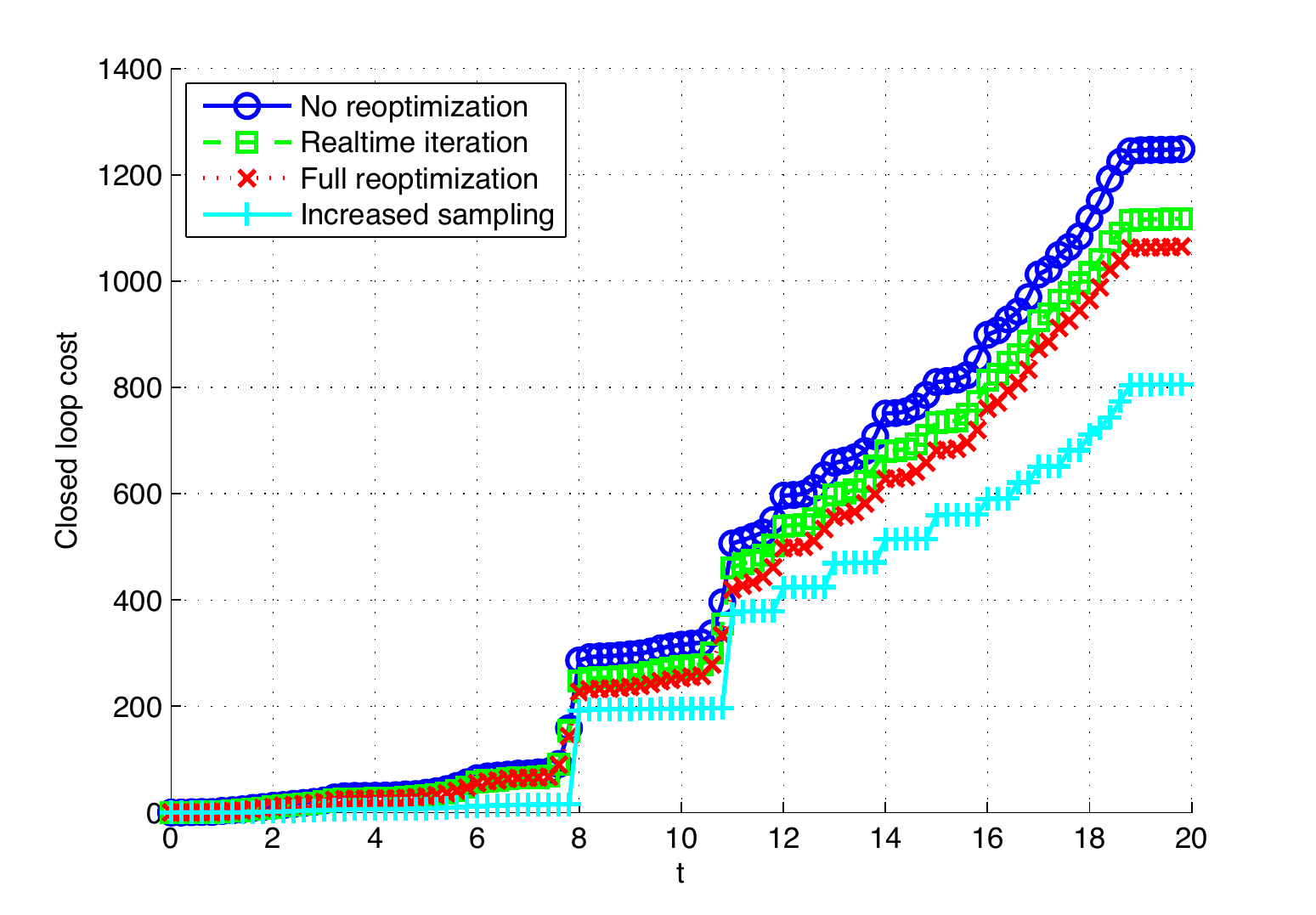}
		\caption{Closed loop costs for different updating techniques}
		\label{Numerical Results:fig:closed loop costs}
	\end{figure}
\end{center}

In addition to measurement errors of the road profile, we then analyzed the impact of additional state deviations due to, e.g., modelling errors, which is uniformly distributed in the interval $[-0.005, 0.005]$, i.e. $\| \statedisturbed \|_{\state} \leq \statedisturbedmax = 0.005$ in Assumption \ref{Setup and Preliminaries:ass:bounded disturbance}, on the closed loop performance regarding different updating techniques. As we suppose that no more precise information is available to our system we utilized the scattered road data as input parameters and calculated the updates due to the state measurements. Since the sensitivity based approach is theoretically unsuitable even without these additional errors, we focused on the realtime iteration and the hierarchical MPC approach. As expected, the updated control laws show an improved performance in terms of the closed loop costs shown in Figure \ref{Numerical Results:fig:closed loop costs}: For the chosen road profile the closed loop costs using realtime iterations with sampling time $T_\control = 0.1$ decrease by approximately $10.5 \%$. While this may appear to be a small improvement given the additional computing costs, we added full reoptimization results which show that at most a reduction of approximately $14.6 \%$ is to be expected for our chosen setting.
\vspace{-0.5cm}
Considering hierarchical MPC, however, allows us to increase the sampling rate which is still covered by our theoretical results. To be comparable, we performed one \glqq D--level\grqq{} step each $0.1s$ and one \glqq A--level\grqq{} step each $0.002s$, i.e. $T_\control = T_\parameter = 0.002s$ and we obtain 50 times as many updates as in the other approaches, cf. \cite{BockDiehlKostinaSchloeder2007} for details on the updating steps. Even without the intermediate \glqq B--level\grqq{} and \glqq C--level\grqq{} steps, we observed an improvement of approximately $35.4 \%$.
\begin{center}
	\begin{figure}[!htb]
		\includegraphics[width=0.47\textwidth]{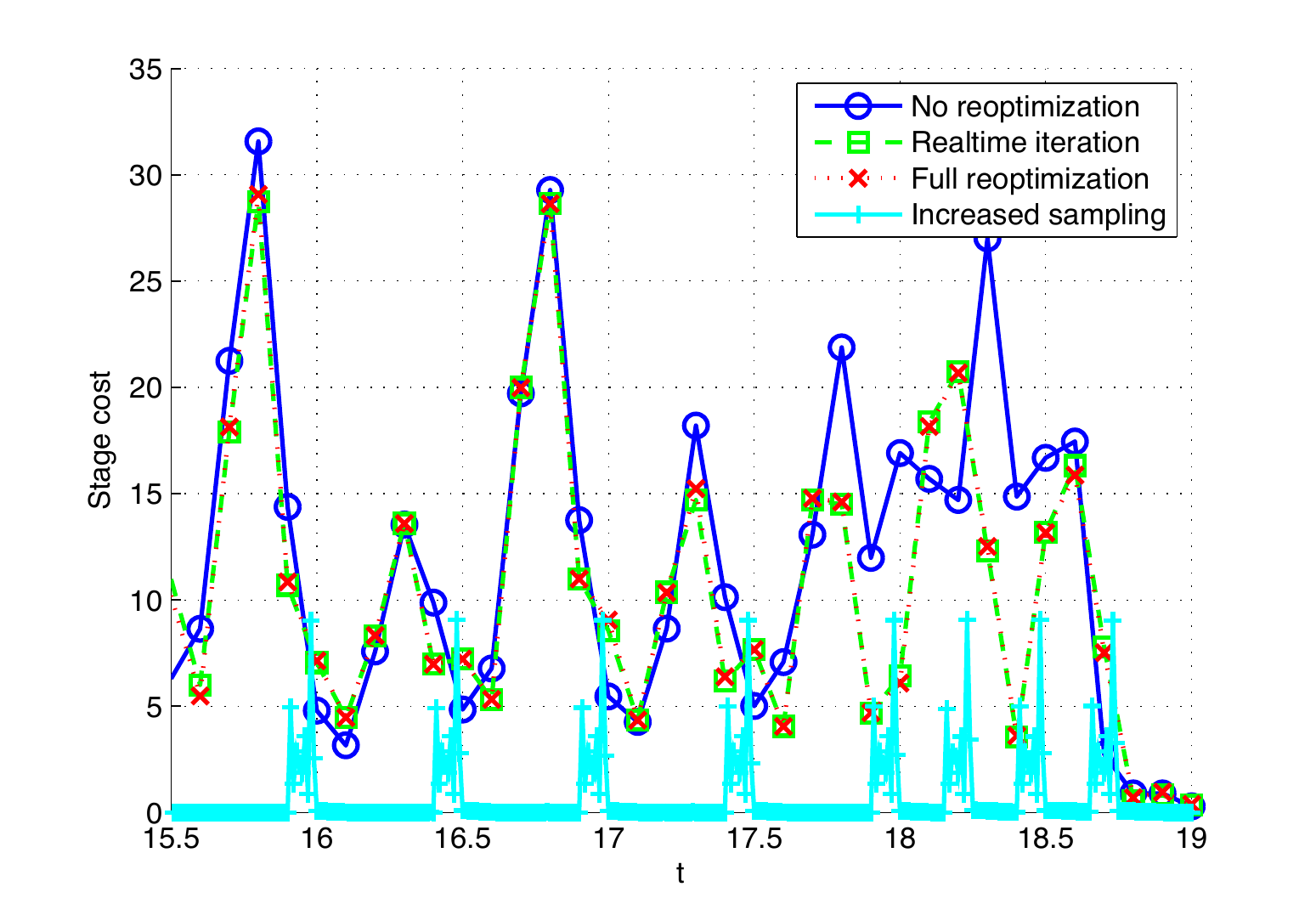}
		\caption{Stage costs for different updating techniques}
		\label{Numerical Results:fig:stage loop costs}
	\end{figure}
\end{center}
\vspace{-0.5cm}
Taking a closer look at the occuring stage costs depicted in Figure \ref{Numerical Results:fig:stage loop costs}, it becomes clear that the faster switching ability allows not only to reduce the costs during jumps of the road but also to recover faster from these jumps, i.e. both the overshoot and the decay rate of excitations are improved.

\section{Conclusions and Outlook}
\label{Section:Conclusions and Outlook}

We presented a feasibility and stability proof for nonlinear model predictive control with abstract updates using relaxed Lyapunov arguments. In particular, we have shown that the known updating techniques via sensitivity information, realtime iterations and hierarchical MPC are covered by our result.\\
Utilizing the jump constants within our assumptions, future research will focus on incorporating model reduction effects on the control into the stability analysis as well as extensions towards hybrid systems by covering integer jumps in variables.

\begin{ack}
This work was partially funded by the German Federal Ministry of Education and Research (BMBF), grant no. 05M10WNA. 
\end{ack}

\bibliographystyle{ifacconf}


\end{document}